\documentclass[a4paper,11pt]{article}
\pdfoutput=1 

\usepackage{jheppub} 

\usepackage[T1]{fontenc} 

\usepackage{pdfsync}
\usepackage{mathptmx}
\usepackage[utf8]{inputenc}
\usepackage[T1]{fontenc}
\usepackage{slashed}
\usepackage{times}
\usepackage{amssymb,amsfonts,amsmath,amsthm}
\usepackage{dsfont,bbm}
\usepackage{dcolumn}
\usepackage{epsf}
\usepackage{graphicx}
\usepackage[caption=false]{subfig}
\usepackage{dsfont}
\usepackage{simplewick}
\usepackage{bm}
\usepackage{eucal}
\usepackage{slashed}
\usepackage[active]{srcltx}
\usepackage[usenames]{color}

\title{\boldmath Complexity of Radon transforms
}

\author[a]{I.~V.~Anikin}
\affiliation[a]{Bogoliubov Laboratory of Theoretical Physics JINR, 141980 Dubna, Russia}

\emailAdd{anikin@theor.jinr.ru}

\abstract{
For the reconstruction problem, the universal representation of inverse Radon transforms
implies the needed complexity of the direct Radon transforms which leads to the additional contributions.
In the standard theory of generalized functions, if the outset (origin) 
function which generates
the Radon image is a pure-real function, as a rule, the complexity of Radon transforms
becomes in question.
In the paper, we discuss the Fourier slice theorem analyzing the degenerated (singular)
points as possible sources of the complexity.
We also demonstrate the different methods 
to generate the needed complexity on the intermediate stage of calculations.
Besides, we show that the introduction of the hybrid (Wigner-like) function
ensures naturally the corresponding complexity.
The discussed complexity provides not only
the additional contribution to the inverse Radon transforms, but also it makes an essential impact
on the reconstruction and optimization procedures within the
frame of the incorrect problems.
The presented methods can be effectively used for the practical tasks of reconstruction problems.
}

\begin{document}
\maketitle
\flushbottom

\section{Introduction}
\label{Intro}

In the modern world, computed tomography technologies (CT-technologies) influence on the many different fields
basically owing to a possibility to investigate the internal composite structure of a object without cutting
and breaking.
From the mathematical viewpoint, CT-technologies are closely associated with the application of both
the direct and inverse Radon transforms \cite{Deans}. It is well-known that due to the inverse Radon transforms
one can visualize the internal structure of the objects under investigation.
The quality of visualization strongly depends on
the inversion procedure of Radon transforms that, as known, is ill-posed and demands
the corresponding regularization, see for example \cite{Anikin:2019oes}.

In \cite{Anikin:2024vto, Anikin:2025lta}, the universal (or unified) representation of inverse Radon transform
has been proposed and studied for any dimension of space.
Thanks for the use of generalized function theory, the mentioned universal representation
involves the new contribution (in comparison with the standard methods based on the Courant-Hilbert identities \cite{Courant-Hilbert})
which definitely effects on the reconstruction and extends the
Tikhonov-like regularization (AC-regularization) used for the solutions of the incorrect task classes.
In particular, within the reconstruction procedure, the additional new contribution
in the unified inverse Radon transforms is given by the integration
with the complex measure that should be compensated, generally speaking, by the complexity of direct Radon image.

In the present paper,  we analyze in detail the Fourier slice theorem which 
plays the extremely important role for the derivation of the universal inverse Radon transforms.
Moreover, we give the interpretations of degenerated (singular) points which require the 
certain regularization and, as a consequence, the corresponding complexity.   

We also propose the different methods
how to generate the needed complexity on the intermediate stage of calculations.
Our methods are based on {\it (i)} the demonstrated non-trivial holonomy
of Radon transforms; 
{\it (ii)} the Fourier series expansion applied within the wide-spread discrete slice approximation \cite{Deans} together with the use of the hybrid (Wigner-like) functions.
The latter are, by definition, the functions of both spacial and momentum coordinates.
For the universal inverse Radon transforms and, then, for the reconstruction procedure,
this trick has been never used up to now.

\section{The Fourier slice theorem and the degenerated (singular) points}
\label{FsTh}

In this section,
we analyse in detail the Fourier slice theorem which is the principal step 
toward the universal inverse Radon transforms.

Let $f(\vec{\bf x})$ with $\vec{\bf x}\in \mathbb{R}^2$ (this is $x$-space) be a original outset
function the internal structure of which is under investigation without
breaking.
For the sake of simplicity, we focus on the two-dimensional coordinate space,
the extension up to an arbitrary dimension is obvious, see for example \cite{Anikin:2025lta}.

\textbf{Theorem:}
{\it the Fourier image of the outset function $f(\vec{\bf x})$ relates to the direct Radon image of the same
function $f(\vec{\bf x})$ through the one-dimensional Fourier transformation with respect to the radial coordinate}:
\begin{eqnarray}
\label{F-t-3-dir}
\boxed{
\mathcal{F}[f](\lambda, \varphi) =\int_{-\infty}^{+\infty}(d\tau) \,e^{-i\lambda\tau}\,
\mathcal{R}[f](\tau,\varphi).
}
\end{eqnarray}

Here and in what follows,
the brackets in the integration measure, see $(d\tau)$, denote that the
corresponding normalization factors have been included and they
are not written explicitly unless it leads to misunderstanding.

On the {\it l.h.s} of (\ref{F-t-3-dir}),
the Fourier image of $f(\vec{\bf x})$ is given by
\begin{eqnarray}
\label{F-t-1-dir}
\mathcal{F}[f](\vec{\bf q})= \int_{-\infty}^{+\infty} d^2 \vec{\bf x} \,
e^{-i\langle\vec{\bf q},\vec{\bf x}\rangle} \,f(\vec{\bf x})
\end{eqnarray}
and
the polar coordinate system:
$\vec{\bf q}\equiv \lambda\, \vec{\bf n}_{\varphi}$ with $|\vec{\bf n}_{\varphi}|=1$
is being implied in $q$-space that is conjugated to $x$-space.
While, on the {\it r.h.s} of (\ref{F-t-3-dir}),
the direct Radon transform of $f(\vec{\bf x})$ is defined as
\begin{eqnarray}
\label{F-t-4-dir}
\mathcal{R}[f](\tau,\varphi)=
\int_{-\infty}^{+\infty} d^2 \vec{\bf x} \, f(\vec{\bf x})
\delta\left( \tau - \langle\vec{\bf n}_{\varphi}, \vec{\bf x}\rangle\right),
\end{eqnarray}
where the argument of delta-function pramaterizes the straightforward line
in $\mathbb{R}^2$.
As seen from (\ref{F-t-3-dir}) the angular dependence of Fourier and Radon images are
coinciding.

To proof the inference of (\ref{F-t-3-dir}), we insert the integral unit into the {\it r.h.s.} of
(\ref{F-t-1-dir}). It reads
\begin{eqnarray}
\label{F-t-1-dir-2}
\mathcal{F}[f](\vec{\bf q})= \int_{-\infty}^{+\infty} d^2 \vec{\bf x} \,
e^{-i\langle\vec{\bf q},\vec{\bf x}\rangle} \,f(\vec{\bf x})
\times
\Big\{
\int_{-\infty}^{+\infty} (dt) \delta\big( t -  \langle\vec{\bf q},\vec{\bf x}\rangle\big)
\Big\}.
\end{eqnarray}
The 
Introducing
\begin{eqnarray}
\label{Par-1}
\xi\stackrel{\text{def.}}{=}\frac{q_2}{q_1}, \quad
z\stackrel{\text{def.}}{=} \frac{t}{q_1},
\end{eqnarray}
one can readily obtain
\begin{eqnarray}
\label{F-t-1-dir-3}
\mathcal{F}[f](q_1, \xi q_1)= \int_{-\infty}^{+\infty} (dz)
e^{-i z q_1}
\Big\{
 \int_{-\infty}^{+\infty} d x_2 \,
f(z-\xi x_2, x_2)
\Big\}
\equiv
\int_{-\infty}^{+\infty} (dz)
e^{-i z q_1}  \, \mathcal{R}[f](z,\xi)
\end{eqnarray}
which can be represented in the form of (\ref{F-t-3-dir}) if the polar system has been used
and the variable $\tau = t/\lambda$
has been introduced instead of $z$, see (\ref{Par-1}).

Notice that the Fourier variable $\lambda\in [0, +\infty]$
(the radial component of $q$-space)
and the Radon radial variable $\tau\in [-\infty, +\infty]$ form the conjugated set.
Hence, the singular point $\lambda=0$ (it is the so-called ``axis'' singularity of the polar system),
leading to the degeneration of angular dependence \cite{Anikin:2025lta},
corresponds to the point $\tau=\pm\infty$.
If the Radon image as a function of arguments is bounded and is finite for some reasons,
the axis singularity is not important.
On the other hand, the point $\lambda=+\infty$ may generate the other type of singularity
and, therefore, the conjugated Radon variable $\tau=0$ should be also regularized
\cite{Anikin:2025lta}.

The mentioned situation appears in the inversion problem where both sides of  (\ref{F-t-3-dir})
have been integrated out over the radial and angular variables.
Indeed, having used the inverse Fourier transform together with (\ref{F-t-3-dir}), we can write down that
\begin{eqnarray}
\label{Inv-F-t-2}
f(\vec{\bf x})&=&
\int_{-\infty}^{+\infty} d^2 \vec{\bf q}\,  e^{+i\langle\vec{\bf q},\vec{\bf x}\rangle} \,\mathcal{F}[f](\vec{\bf q})
\Big|_{\vec{\bf q}=\lambda\vec{\bf n}_{\varphi}}
\\
\label{Inv-F-t-2-2}
&=&
\int_{0}^{+\infty} d\lambda \lambda \int_{\text{f. r.}} d \varphi\,
e^{+i\lambda\langle \vec{\bf n}_{\varphi}, \vec{\bf x}\rangle}\,
\mathcal{F}[f](\lambda, \varphi),
\end{eqnarray}
where $\text{``f. r.''}$ implies that the angular integration measure covers the full regions of variations.
In the inverse Fourier transform, see (\ref{Inv-F-t-2}), the integration in $q$-space should be 
understood in a sense of the principle value integration. It is needed to avoid the problems with the
existence of improper integrations (see, for example, \cite{Gelfand:1964}).  

Hence, the  use of the Fourier slice theorem, see (\ref{F-t-3-dir}), yields
\footnote{$\epsilon$ as a subscript of $f$ denotes the $\epsilon$-regularization that
should be used in (\ref{Inv-F-t-2-3}).}
\begin{eqnarray}
\label{Inv-F-t-2-3}
f_\epsilon(\vec{\bf x}) =
\int_{\text{f. r.}} d\varphi\,
\int_{-\infty}^{+\infty} (d\eta)\,
\mathcal{R}[f](\eta + \langle \vec{\bf n}_{\varphi}, \vec{\bf x}\rangle, \varphi)
\times
\Big\{\int_{0}^{+\infty} d\lambda \lambda \,
e^{-i\lambda\eta}\Big|_{\epsilon\text{-reg.}}
\Big\},
\end{eqnarray}
where, in the factorized $\lambda$-integration,  $``\epsilon\text{-reg.}"$ denotes the necessary regularization of $\lambda$-integration. 

Since, in (\ref{Inv-F-t-2-3}), the $\lambda$-integration has been factorized out 
it can be independently regularized as the regularization of $\eta$-pole
(which occurs after the integration over $\lambda$):
\begin{eqnarray}
\label{Pole-reg}
\delta_+(\eta)=\int_{0}^{+\infty} (d\lambda)  \, e^{-i\lambda\eta}\Big|_{\epsilon\text{-reg.}}
\equiv
\int_{0}^{+\infty} (d\lambda)  \, e^{-i\lambda(\eta - i\epsilon)}=
\frac{1}{\eta - i\varepsilon}.
\end{eqnarray}
Alternatively, for regularization, 
we can make a replacement: $\eta\to \eta - i\varepsilon$ in the integration over $\eta$
in (\ref{Inv-F-t-2-3}). It results in the complex continuation of the Radon radial component
the necessity of which agrees with the dual Radon transforms, 
see \cite{Anikin:2019oes, Anikin:2024vto, Anikin:2025lta}.

To conclude this section, the imaginary part of $\delta_+(\eta)$ of (\ref{Pole-reg}) 
disappears provided
the full region of angular integration is considered
 \cite{Anikin:2019oes, Anikin:2024vto, Anikin:2025lta}.
It is worth to stress that the full regions of angular and radial integrations in the polar system
correspond to the full region of integration assumed, from the very beginning,
in the Cartesian system of $q$-space, see (\ref{Inv-F-t-2}).
Moreover, it is already clear that (\ref{Inv-F-t-2-3}) establishes a basis for the universal inversion of Radon transforms
in the case of the restricted angular integration \cite{Anikin:2025lta}.

\section{The linearity and outset functions with defects}
\label{Non-lin}

The linearity of the direct Radon transforms:
\begin{eqnarray}
\label{Lin-1}
\mathcal{R}[f_1+f_2](\tau, \varphi) = 
\mathcal{R}[f_1](\tau, \varphi) + \mathcal{R}[f_2](\tau, \varphi),
\end{eqnarray}
is one of basic properties that should be satisfied. 
The condition (\ref{Lin-1}) is valid
provided, in particular, both outset functions have the full-region supports
and they have the regular properties. 
Besides, the corresponding integration measures that define the transformations, as a rule,
cover also the full region \cite{Deans, Gelfand:1964, GGV}. 

Based on the applications of Radon transforms in CT-technology, 
the outset function can be presented as a combination of functions
which are differently localized in $x$-space.
In this section, we dwell on this case, 
dictated by practice. 
The latter can also lead to the complexity of Radon transforms.  
 
At the beginning, 
let $f_1(\vec{\bf x})$ and $f_2(\vec{\bf x})$ be 
$g_1(\vec{\bf x}) \Theta_I(\vec{\bf x})$ and $g_2(\vec{\bf x}) \Theta_{III}(\vec{\bf x})$, respectively, 
where $\Theta_{I, III}$ are the characteristic indicators defined as
\begin{eqnarray}
\label{Char-Theta} 
\Theta_{I}(\vec{\bf x}) =
\left\{
  \begin{aligned}
  1,\quad &\text{if} \quad \vec{\bf x} \in \{\Omega_I | x_i \geq 0 \}\\
  0, \quad &\text{if} \quad \vec{\bf x} \notin \{\Omega_I | x_i \geq 0 \}
  \end{aligned}
\right., \quad 
\Theta_{III}(\vec{\bf x}) =
\left\{
  \begin{aligned}
  1,\quad &\text{if} \quad \vec{\bf x} \in \{\Omega_{III} | x_i < 0 \}\\
  0, \quad &\text{if} \quad \vec{\bf x} \notin \{\Omega_{III} | x_i < 0 \}
  \end{aligned}
\right.\,.
\end{eqnarray}
As mentioned, in the direct Radon transform (\ref{F-t-4-dir}), the argument of delta-function gives the line parametrization 
as a function of the radial $\tau$ and angular $\varphi$ variables which define the position of the 
straightforward line in $x$-space. On the other hand, the direct Radon transform can be reduced to the 
curve-linear integral of the first kind over the parameter $s$ that determines the integration measure 
along the given line, see \cite{Anikin:2025lta}.

If the outset function is combined as 
\begin{eqnarray}
\label{g1-g2-com}
G(\vec{\bf x})=g_1(\vec{\bf x}) \Theta_I(\vec{\bf x}) + g_2(\vec{\bf x}) \Theta_{III}(\vec{\bf x}), 
\end{eqnarray}
the functional associated with the Radon transform is given by 
\begin{eqnarray}
\label{F-RT-1}
\mathcal{R}[G](\tau,\varphi)\stackrel{\text{def.}}{=}
\int_{-\infty}^{+\infty} d^2 \vec{\bf x} \, 
\big\{
g_1(\vec{\bf x}) \Theta_I(\vec{\bf x}) + g_2(\vec{\bf x}) \Theta_{III}(\vec{\bf x}) 
\big\}
\delta\left( \tau - \langle\vec{\bf n}_{\varphi}, \vec{\bf x}\rangle\right).
\end{eqnarray}
Hence, one can see that the components of $\vec{\bf x}$ should 
semultaneously 
satisfy two conditions:
{\it (a)} the delta-function singles out the points which lie on the given straightforward line
and {\it (b)} the points have to belong either $\Omega_I$ or $\Omega_{III}$ 
due to the characteristic indicators.

Focusing on the case of $\tau > 0$ and $\varphi \in [0,\, \pi/2]$, we get the following
\footnote{In the similar manner we can treat the case of $\tau > 0$ and $\varphi \in [\pi,\, 3\pi/2]$.}
\begin{eqnarray}
\label{F-RT-2}
\mathcal{R}[G](\tau,\varphi)\Big|^{\tau > 0}_{\varphi \in [0,\, \pi/2]}=
\int_{-\infty}^{+\infty} d^2 \vec{\bf x} \,
g_1(\vec{\bf x}) \Theta_I(\vec{\bf x}) 
\,
\delta\left( \tau - \langle\vec{\bf n}_{\varphi}, \vec{\bf x}\rangle\right)
=\mathcal{R}[g_1](\tau,\varphi).
\end{eqnarray}
Then, we make a shift: $\varphi\to \varphi + \pi$ in $\mathcal{R}[G](\tau,\varphi)$ of 
(\ref{F-RT-1}), we get
\begin{eqnarray}
\label{F-RT-3}
\mathcal{R}[G](\tau,\varphi + \pi)\Big|^{\tau > 0}_{\varphi \in [0,\, \pi/2]}=
\int_{-\infty}^{+\infty} d^2 \vec{\bf x} \,
g_2(\vec{\bf x}) \Theta_{III}(\vec{\bf x}) 
\,
\delta\left( \tau - \langle\vec{\bf n}_{\varphi+\pi}, \vec{\bf x}\rangle\right)
=\mathcal{R}[g_2](\tau,\varphi+\pi).
\end{eqnarray}
As the next step, we perform the other shift on $\pi$ in $\mathcal{R}[g_2]$ 
on the {\it r.h.s.} of (\ref{F-RT-3}) to obtain
\begin{eqnarray}
\label{F-RT-4}
\int_{-\infty}^{+\infty} d^2 \vec{\bf x} \,
g_2(\vec{\bf x}) \Theta_{III}(\vec{\bf x}) 
\,
\delta\left( \tau - \langle\vec{\bf n}_{\varphi+\pi+\pi}, \vec{\bf x}\rangle\right)
=\mathcal{R}[g_2](\tau,\varphi+\pi+\pi)=0
\end{eqnarray}
owing to the characteristic indicator $\Theta_{III}$.
However, on the other hand, the shift on $2\pi$ of $\mathcal{R}[G](\tau,\varphi)$, see 
(\ref{F-RT-1}), leads to 
\begin{eqnarray}
\label{F-RT-5}
\mathcal{R}[G](\tau,\varphi+2\pi)\Big|^{\tau > 0}_{\varphi \in [0,\, \pi/2]}
=\mathcal{R}[g_1](\tau,\varphi)
\end{eqnarray}
and, therefore, we have 
\begin{eqnarray}
\label{F-RT-6}
\boxed{
\mathcal{R}[g_2](\tau,\varphi+\pi+\pi)
\not=\mathcal{R}[g_1](\tau,\varphi).
}
\end{eqnarray}
The condition (\ref{F-RT-6}) points to the existence of the non-trivial holonomy
which can generate the complexity discussed in the present paper.
In addition, owing to the different properties of $g_1$ and $g_2$ of (\ref{g1-g2-com}),
it is not excluded that 
\begin{eqnarray}
\label{right-left-lim}
\lim_{\tau\to 0-} \mathcal{R}[g_2](\tau, \varphi) \not= 
\lim_{\tau\to 0+} \mathcal{R}[g_1](\tau, \varphi),
\end{eqnarray}
giving the disconnection 
which
is essential for the  contributions of $F_A, \widetilde{F}_A$, see 
below (\ref{Inv-R-F-1}).

Having used the mentioned demonstration, we are in a position to discuss 
the outset function which contains some ``defect''.
We assume that 
\begin{eqnarray}
\label{g1-g2-com-W1}
\tilde G(\vec{\bf x})=
\big[ g_1(\vec{\bf x}) + g_2(\vec{\bf x})\big]\Theta_I(\vec{\bf x}) 
+ g_2(\vec{\bf x}) \Theta_{III}(\vec{\bf x}), 
\end{eqnarray}
where $g_1$ plays now a role of some ``defect'' that is localized in $\Omega_I$.

We again consider the rotation on $2\pi$ in the direct Radon transforms by 
two ways:
{\it (i)} $\varphi\to \varphi + 2\pi$ and $\varphi\to\varphi + \pi +\pi$ (with the 
calculation of Radon functional at the intermediate point $\varphi+\pi$).
Having used the results of (\ref{F-RT-2})-(\ref{F-RT-6}), we readily derive that  
\begin{eqnarray}
\label{F-RT-7}
\boxed{
\mathcal{R}[g_2](\tau,\varphi+\pi+\pi)
\not=\mathcal{R}[g_1+g_2](\tau,\varphi).
}
\end{eqnarray}
This allows us to extract the influence of the ``defect'' describing by $g_1$ on the 
common ``background'' in $\Omega_I$ and $\Omega_{III}$ given by $g_2$. 

To summarize this section, we have presented the evidences for the non-trivial holonomy 
of Radon transforms. In particular, we have suggested the way how to single out 
the effects from the so-called ``defects'' hidden into the given functions
describing the object under CT-investigation.

\section{Two-dimensional space versus three-dimensional space}
\label{Dim-2-3}

In \cite{Anikin:2025lta}, the universal and unified representation for the inversion of Radon transforms has been presented
where the Courant-Hilbert identities \cite{Courant-Hilbert} have not been applied.  We remind that the Courant-Hilbert identities are based on the
use of the Green formulae the different forms of which depend on either even or odd dimension of space \cite{Gelfand:1964, GGV}.
For the practical uses (for example, in medicine), the most needed cases correspond to $\mathbb{R}^3$ and $\mathbb{R}^2$ spaces.

However, as dictated by the application of CT-technology, the different algorithms have been
designed for the two-dimensional reconstructions which work with the
corresponding transverse projections of the three-dimensional object.

In other words, if the object, which is needed to be reconstructed, is described by the three-dimensional
origin outset function $f(x_1, x_2, x_3)$, we perform the finite number of transverse sections, say, regarding $x_3$-axis:
\begin{eqnarray}
\label{S-1}
f(x_1, x_2, x_3) \Longrightarrow \Big\{  f(x_1, x_2, x_{3_1}), \, f(x_1, x_2, x_{3_2}),\, ....,\, f(x_1, x_2, x_{3_n}) \Big\},
\end{eqnarray}
where the number of sections should be defined by the experiment.
The segmentation given by (\ref{S-1}) is known as  {\it a discrete slice method} known in CT-technology.
In this case, we always deal with the direct and inverse Radon transforms determined on the two-dimensional space formed by $(x_1, x_2)$
while the third discrete coordinate plays a role of the external parameter \cite{Deans}.

\section{The hybrid Wigner-like function and the complex Radon transforms}
\label{Complex-R}

As explained in a series of papers \cite{Anikin:2019oes, Anikin:2024vto, Anikin:2025lta},
the regularized inverse Radon operator involves two contributions. If one of them is related to the
real integration measure, another is associated with the imaginary integration measure
(that is a result of the Cauchy theorem).

According the scenario described in \cite{Anikin:2024vto} the term of the inverse Radon operator with
the imaginary measure is extremely
important for the opimization procedure which works with the condition inspired by the corresponding norms.
On the other hand, in QFT, even the direct Radon transform, which is linked to the corresponding
transverse-momentum dependent parton distributions, can possess the imaginary part owing to the interactions
in the correlators \cite{Anikin:2019oes}.
In this section, excepting QFT-models, we study the natural source of complexity which appears in the Radon transforms.

In this section, we now begin with $\mathbb{R}^3$-space where 
the coordinate system has been defined.
Let $f(\vec{\bf x})$, with $\vec{\bf x}\in \mathbb{R}^3$, be a outset (origin) function which is usually
well-localized. In the reconstruction problem, the form of $f(\vec{\bf x})$ has to be restored
owing to the inverse Radon transform.

Since the integral representation of universal inverse Radon transform involves the additional term which should be related to the
complex integrand given by the direct Radon transform,
we first introduce the hybrid (Wigner-like) function
\footnote{In the particular case of QFT, the Wigner function corresponds to the quasi-probability distribution which is
actually a real function only.}
associated with the outset function.
Namely, we
perform the direct Fourier transform regarding one of three coordinates 
applied to $f(\vec{\bf x})$.
We have
\begin{eqnarray}
\label{D-F-1}
F(x_1, x_2; k)= \int\limits_{-\infty}^{+\infty} (dx_3)\, e^{-ikx_3} \,f(x_1, x_2, x_3),
\end{eqnarray}
where, as above-used, the irrelevant normalization has been absorbed into the integral measure.
In (\ref{D-F-1}), it is clear that $f(x_1, x_2, x_3)\in \Re{\rm e}$ by definition, while $F(x_1, x_2; k)\in \mathbb{C}$ already.
Hence, in the discrete slice method, see (\ref{S-1}), we have the following
\begin{eqnarray}
\label{D-F-2}
\widetilde{F}(x_1, x_2; k)= \sum_n\, e^{-ik x_{3_n}} \,f(x_1, x_2, x_{3_n}), \quad \widetilde{F}(x_1, x_2; k) \in \mathbb{C},
\end{eqnarray}
which is nothing but the Fourier series expansion.

Then, we calculate the Radon image of the Fourier $F$-functions as
\begin{eqnarray}
\label{D-R-F-1}
\mathcal{R}[F; \widetilde{F}](\tau,\varphi; k)=
\int_{-\infty}^{+\infty} d^2 \vec{\bf x} \,
  \left[ \begin{array}{c}
F(x_1, x_2; k)\\
\widetilde{F}(x_1, x_2; k)
  \end{array}\right]
\,
\delta\left( \tau - \langle\vec{\bf n}_{\varphi}, \vec{\bf x}\rangle\right),
\end{eqnarray}
where $\mathcal{R}[F; \widetilde{F}](\tau,\varphi; k) \in \mathbb{C}$ and
the unit vector is given by $\vec{\bf n}_{\varphi}=(\cos\varphi, \sin\varphi)$.
Notice that in (\ref{D-R-F-1}) the integration measure corresponds to the two-dimensional space while
the momentum $k$ plays a role of the external parameter.

Hence, the inverse Radon transform in the universal form is given by
\begin{eqnarray}
\label{Inv-R-F-1}
\boxed{
F_\epsilon(x_1, x_2; k) = F_S(x_1, x_2; k) + F_A(x_1, x_2; k),
}
\end{eqnarray}
where
\begin{eqnarray}
\label{Inv-F-t-2-4-S2}
&&F_S(x_1, x_2; k) =
-\int_{-\infty}^{+\infty} (d\eta)\, \frac{\mathcal{P}}{\eta^2}\,
\int d\mu(\varphi) \, \mathcal{R}[F](\eta + \langle \vec{\bf n}_\varphi, \vec{\bf x}\rangle, \varphi; k)
\end{eqnarray}
and
\begin{eqnarray}
\label{Inv-F-t-2-4-A2}
&&F_A(x_1, x_2; k) =
-i\pi \int d\mu(\varphi) \,  \frac{\partial}{\partial\eta} \,
\mathcal{R}[F](\eta + \langle \vec{\bf n}_\varphi, \vec{\bf x}\rangle, \varphi; k)\Big|_{\eta=0}.
\end{eqnarray}
The similar expressions can be written for $\widetilde{F}$-function due to the trivial replacement.

In practical applications, the values of $\mathcal{R}[F; \widetilde{F}]$ are usually known from the experiment (or from the observation).
Restoring the function $F$ or $\widetilde{F}$ from $\mathcal{R}[F; \widetilde{F}]$, one can then reconstruct the function $f$. 
In this connection, based on the slice method (\ref{S-1}) the introduction of 
$\widetilde{F}$ aggregates all two-dimensional slices with the help of the continuous 
parameter $k$, see (\ref{D-F-2}) which is a very convenient parameter for the 
optimization procedure \cite{Anikin:2024vto}.
  
Therefore, for the reconstruction, the ultimate scheme can expressed as
\begin{eqnarray}
\label{F-1}
\boxed{
\mathcal{R}[F; \widetilde{F}](\tau, \varphi; k) \stackrel{{\cal R}^{-1}}{\longrightarrow}
\{ F(x_1, x_2; k), \widetilde{F}(x_1, x_2; k) \}
 \Longrightarrow f(x_1, x_2, x_3).
}
\end{eqnarray}
This scheme illustrates the principle result which shows us that the necessary step to
obtain the complexity of direct Radon transform is given thanks for the introduction and
the transition to the hybrid functions $F(x_1, x_2; k)$ or $\widetilde{F}(x_1, x_2; k)$.
Notice that the use of Fourier images in the forms of the hybrid $F$-functions, see (\ref{D-F-1}) and (\ref{D-F-2}),
gives one of possibilities to deal with the complex direct Radon transforms.
We remind that $F_S$ and $\widetilde{F}_S$ are determined by $\Re{\rm e}\big\{\mathcal{R}[F; \widetilde{F}]\big\}$
and $F_A$ and $\widetilde{F}_A$ are related to $\Im{\rm m}\big\{\mathcal{R}[F; \widetilde{F}]\big\}$
garanteeing the new important contribution within the reconstruction problem.

\section{Conclusions}
\label{Cons}

In the paper, we have described the methods where
the complexity of the Radon transform appears naturally.
In the connection with the reconstruction problem, the complex Radon transforms
are started to be extremely important ones. This is because of the additional contribution existence
that plays a crucial role in the improved visualization procedure \cite{Anikin:2024vto}.
 
We have demonstrated the conditions for the non-trivial holonomy of Radon transforms 
leading to the complexity of transformations.
In addition, we have also presented the method to extract the ``defect'' hidden into
the common ``background'' function. This opens a window for the 
development of new reconstruction schemes.  

Furthermore, we have presented 
the method which is practically identical to the Fourier series expansion applied
within the frame of the discrete two-dimensional slice approximation \cite{Deans}.
In its turn, the discrete two-dimensional slice approximation is a wide-spread approximation in the
medical CT-technologies.
The proposed method is based on the use of the hybrid (Wigner-like) functions
depending simultaneously on the spacial and momentum coordinates.
The introduced hybrid function as a complex function can be considered as a 
very convenient tool for the universal inverse Raton transforms.
This trick has been first applied in the context of the inversion procedure of Radon transforms.

\section*{Acknowledgements}
We thank A.I.~Anikina, V.A.~Osipov and O.I.~Streltsova for useful and illuminating discussions.
The special thanks go to the colleagues from the South China Normal University (Guangzhou) for the useful discussions and a very warm hospitality.

\section*{Conflict of interest}

The author declares no conflict of interest.

\section*{Data Availability Statement}
This manuscript has no associated data or the data will not be deposited.


\end{document}